\documentclass{ifacconf}
\usepackage{savesym}
\savesymbol{AND}

\usepackage{algorithm}
\usepackage{algorithmic}

\usepackage{natbib}
\usepackage{mathrsfs}
\usepackage{amsmath}
\usepackage{amssymb}

\newcommand{\IM}{\mathrm{Im}\mbox{ } }

\newcommand{\Rank}{\mathrm{rank}\mbox{ } }
\newcommand{\SPAN}{\mathrm{Span}}
\newtheorem{Definition}{Definition}
\newtheorem{Theorem}{Theorem}
\newtheorem{Lemma}{Lemma}

\newtheorem{Remark}{Remark}

\newtheorem{Example}{Example}

\newtheorem{definition}{Definition}

\usepackage{color}

\usepackage{graphicx}
\usepackage{epsfig}
\usepackage{color}
\usepackage{subfigure}
\usepackage{epic}
\usepackage{eepic}

\usepackage{pdfpages}

\begin{document}


\newcommand{\diag}{\mathrm{diag}}

\begin{frontmatter}

\title{Moment matching for bilinear systems with nice selections}
\author[FIRST]{Mih\'aly Petreczky}
\author[SECOND]{Rafael Wisniewski and John Leth}

\address[FIRST]{
Centre de Recherche en Informatique, Signal et Automatique de Lille (CRIStAL) {\tt\small mihaly.petreczky@ec-lille.fr}
}
\address[SECOND]{
Aalborg University, Dept. of Electronic Systems,
Fr. Bajers Vej 7, C3-211, DK-9220 Aalborg {\O}st, Denmark
\texttt{\{raf,jjl\}@es.aau.dk}
}

\begin{abstract}
The paper develops a method for model reduction of bilinear control systems. It leans upon the observation that the input-output map of a bilinear system has a particularly simple Fliess series expansion. Subsequently, a model reduction algorithm is formulated such that the coefficients of Fliess series expansion for the original and reduced systems match up to certain predefined sets - nice selections. Algorithms for computing matrix representations of unobservability and reachability spaces complying with a nice selection are provided. Subsequently, they are used for calculating a partial realization of a given input-output map.
\end{abstract}

\end{frontmatter}

\section{Introduction}
 The paper deals with model reduction of bilinear systems. Bilinear systems appear in many applications, and there is an extensive literature on their
 analysis and control \cite{Isi:Nonlin,ElliotBook}. However, bilinear models occurring in practice are often of such large dimensions that application of existing control synthesis and analysis methods is not
 feasible. For this reason, it is of interest to investigate model reduction algorithms for bilinear systems.

\textbf{Problem Formulation}

In the sequel we give a brief description of the system to be studied, for details we refer to Section~\ref{sec:pre}. Consider a bilinear system of the form
\begin{equation*}
\dot x = A_0x+\sum_{i=1}^{m} (A_ix)u_i, \quad y = Cx, 
\end{equation*}
with $A_0,\ldots,A_m \in \mathbb{R}^{n \times n}$, $C \in \mathbb{R}^{p \times n}$, and where $x$
is the state trajectory, $u = (u_1, \hdots, u_m)$ is the input trajectory and $y$ is the output trajectory. 
For a fixed initial state $x_0$, consider the input-output map $\mathcal{Y}_{x_0}$ of the system 
induced by the initial state $x_0$. Specifically, $\mathcal{Y}_{x_0}(u,t)$ is the output $y(t)$ of the system 
at time $t$, corresponding to the input $u$ and initial state $x(0)=x_0$.

In this paper, we propose an algorithm for finding a bilinear system of reduced order:
\begin{equation*}
 \dot {\bar{x}} = \bar{A}_0\bar{x}+\sum_{i=1}^{m} (\bar{A}_i\bar{x})u_i, \quad \bar{y} = \bar{C}\bar{x} 
\end{equation*}
with $\bar{A}_0,\ldots,\bar{A}_m \in \mathbb{R}^{\bar{n} \times \bar{n}}$ and $\bar{C} \in \mathbb{R}^{p \times \bar{n}}$,  $\bar{n} < n$,
and an initial state $\bar{x}_0$ such that the input-output map $\mathcal{\bar{Y}}_{\bar{x}_0}$ of 
the reduced order system 
induced by the initial state $\bar{x}_0$ is close to $\mathcal{Y}_{x_0}$. 
More precisely, the algorithm guarantees that:
\begin{itemize}
\item{\textbf{(i)}}
 $\mathcal{Y}_{x_0}(u,t)=\mathcal{\bar{Y}}_{\bar{x}_0}(u,t)$
for all $u \in \mathcal{U}$ and $t \in [0,T]$, where $\mathcal{U}$ is a priori chosen subset of input
signals, and 
 \item{\textbf{(ii)}} The Euclidian 2-norm
  $\|\mathcal{Y}_{x_0}(u,t) - \mathcal{\bar{Y}}_{\bar{x}_0}(u,t)\|$ is small for small $t$ and $u$.
\end{itemize}
In the paper, we will provide a concrete error bound for $\|\mathcal{Y}_{x_0}(u,t) - \mathcal{\bar{Y}}_{\bar{x}_0}(u,t)\|$ as a function of the magnitude of $t$ and $u$. Moreover, we will discuss how to
choose the class of inputs $\mathcal{U}$ in item \textbf{(i)} above. 
The algorithm resembles that in \cite{MertTAC} developed for the model reduction of linear switched systems. In a nutshall, it relies on \emph{matching} a number of coefficients of the Fliess series (or functional) expansion of
 $\mathcal{Y}_{x_0}$ and $\mathcal{\bar{Y}}_{\bar{x}_0}$, i.e. $\bar{\Sigma}$ is constructed in such a way that a 
certain coefficients of the Fliess-series expansion of $\mathcal{Y}_{x_0}$ and $\mathcal{\bar{Y}}_{\bar{x}_0}$ are the same.
The set of indices of these coefficients  is called a \emph{nice selection}, and it is a subset of sequences of elements from $\{0,1,\ldots,m\}$. 
 Intuitively, various choices of nice selections correspond to choosing those Fliess-series coefficients
which contribute the most to responses to specific input trajectories. In particular, for certain inputs, the coefficients specified by nice selections completely determine the 
output response of the system. In short, we use the nice selection to keep the dimension of the reduced system low by pinpointing our attention on specific input trajectories.  

\textbf{Prior work}
 To the best of our knowledge, the contribution of the paper is new. 
  Model reduction of bilinear systems is an established topic, without claiming completeness, we mention \cite{Bai2006406,Zhang2002205,Wang20121231,Xu20154467,BreitenDamm,BennerSiam,FlagGugercin,BilinearMomentMatching1,BilinearMomentMatching2,BilinearMomentMatching5}. In particular,
 moment matching methods for bilinear systems were proposed in \cite{Bai2006406,BilinearMomentMatching1,BilinearMomentMatching2,BilinearMomentMatching5,BreitenDamm,BennerSiam,FlagGugercin,Wang20121231} and for general non-linear systems in \cite{Astolfi10}. The current paper represents another version of moment matching for bilinear systems.
In \cite{BilinearMomentMatching5,BreitenDamm,BennerSiam} the concept of moment matching at some frequencies $\sigma_1,\ldots,\sigma_k$ was defined. The papers \cite{Bai2006406,BilinearMomentMatching1,BilinearMomentMatching2} correspond to moment matching at zero frequency ($\sigma=0$).
The current paper proposes the notion of $\gamma$-partial realization for some nice selection $\gamma$, which is a generalization of moment matching 
at $\infty$ proposed in \cite{BreitenDamm}. The precise relationship is explained in Remark \ref{comp}. 
Since the notion of $\gamma$-partial realization is more general than moment matching at $\infty$, 
the algorithm of this paper is rather different from those of \cite{Bai2006406,BilinearMomentMatching1,BilinearMomentMatching2,BilinearMomentMatching5,BreitenDamm,BennerSiam,FlagGugercin,Wang20121231}. In particular, we cannot reduce the problem to computing classical Krylov-subspaces, and hence we cannot use 
the corresponding rich mathematical structure. As a result, the numerical 
issues of the proposed algorithm are much less clear than those of the algorithms cited above.
Moreover,the interpretation of the model reduction procedure as an approximation in $H_2$ does not carry over to the framework of the current paper either.
The relationship between the proposed method and moment matching at other frequencies remains the  topic of future research. The same holds for the relative merits of \cite{Bai2006406,BilinearMomentMatching1,BilinearMomentMatching2,BilinearMomentMatching5,BreitenDamm,BennerSiam,FlagGugercin,Wang20121231}. At this stage, it is not clear in which situations the approach proposed in this paper is more usefuk than the algorithms from \cite{Bai2006406,BilinearMomentMatching1,BilinearMomentMatching2,BilinearMomentMatching5,BreitenDamm,BennerSiam,FlagGugercin,Wang20121231}. 
An advantage of the proposed method is that it allows to choose the order of the reduced-order system a-priori, by choosing the number of elements of a nice selection, see
Remark \ref{choice:order} of  Section~\ref{Sec:Algorithms}. 
In addition, to the best of our knowledge, this paper is the first to present a characterization of
those inputs, for which the corresponding input-output behavior is preserved by moment matching. In this sense, the paper follows the spirit of \cite{Astolfi10}, although the technical details of the definition are quite different.
Furthermore,
the paper provides some error bounds, which were absent from the existing literature. 

\textbf{Contents of the paper}
This work is organized as follows. In Section~\ref{sec:pre}, we introduce the notation used throughout the paper; subsequently, we recall the highlights of the realization theory for bilinear systems. In Section~\ref{sect:nice:def}, we introduce the concepts of a column nice selection and a row nice selection, which we use for formulating a convenient partial realization of bilinear systems. Subsequently, we characterize the level at which a bilinear system - a partial realization of $f$ - approximates $f$. Lastly, in Section~\ref{Sec:Algorithms}, we equip the reader with algorithms for computing partial realizations for nice selections.

\section{Preliminaries}\label{sec:pre}
\subsection{Notation}

Denote by $\mathbb{N}$ the set of natural numbers including $0$, and by $\mathbb{R}_+$ the set $[0,+\infty)$ of nonnegative real numbers. The symbol $\|\cdot\|$ will denote the Euclidean 2-norm when applied to vectors and induced 2-norm when applied to matrices. If $F$ is a function between function spaces we write $F[u]$ in place of $F(u)$ to distinguish the arguments e.g., $F[u](t)$ indicate that $F[u]$ is function of $t$.
%
%
%

In addition, $AC(\mathbb{R}_+, \mathbb R^n)$ denotes the set of \emph{absolutely continuous maps}, and $L_{loc}(\mathbb{R}_+, \mathbb R^n)$ the set of \emph{Lebesgue measurable maps} which are integrable on any compact interval. The time derivative of $x\in AC(\mathbb{R}_+, \mathbb R^n)$ is denoted $\dot{x}$ where it is implicitly understood that this notation indicates almost everywhere differentiable.
%

Let $Q$ denote 
the set $\{0,1,\dots,m\}$, and $Q^*$ the set of finite sequences of elements of $Q$ together with the empty sequence $\epsilon$. 
Let $w=q_1q_2 \cdots q_k \in Q^*$ with $q_1,\dots,q_k \in Q$, $k>0$ and $A_{q_i} \in \mathbb{R}^{n \times n}$, $i=1,\dots,k$. Then the matrix $A_w$ is defined as
	\begin{equation*} 
		A_w=A_{q_k}A_{q_{k-1}}\cdots A_{q_1}.
	\end{equation*}
By convention, if $w= \varepsilon$, then $A_\varepsilon$ is the identity matrix.
A bilinear system is a (control) system of the form
	\begin{subequations}\label{eq:BSSform}
		\begin{align} 
			\dot{x}(t)&=A_{0}x(t)+\sum_{i=1}^{m} (A_i x(t)) u_i(t), ~ x(0)=x_0 \label{sdyn} \\
			y(t)&=Cx(t)\label{out} 
		\end{align}
	\end{subequations}
with $A_i \in \mathbb{R}^{n \times n},~i\in Q$, 
$C \in \mathbb{R}^{p \times n}$, and where $u=(u_1,\ldots,u_m)  \in L_{loc}(\mathbb{R}_+, \mathbb R^m)$ is the input trajectory, $x \in AC(\mathbb{R}_+, \mathbb R^n)$ is the state trajectory, and $y \in AC(\mathbb{R}_+, \mathbb R^p)$ is the output trajectory. 

The notation
	\( 	\Sigma=(p,m,n, \{A_i\}_{i\in Q}, C ,x_0)  \)
	or simply $\Sigma$, is used as short-hand representations for a bilinear system of the form \eqref{eq:BSSform}. The number $n$ is the \emph{dimension (or order)} of $\Sigma$ and is sometimes denoted by $\dim \Sigma$.

	The \emph{input-to-state} map $X_{\Sigma,x}$ and \emph{input-to-output} map $Y_{\Sigma,z}$ of $\Sigma$ are the maps 
	\begin{align*}
		X_{\Sigma,z}: L_{loc}(\mathbb{R}_+,\mathbb{R}^m) & \rightarrow AC(\mathbb{R}_+,\mathbb{R}^n); \quad
		u & \mapsto &~ X_{\Sigma,z}[u], \\ 
		Y_{\Sigma,z}: L_{loc}(\mathbb{R}_+,\mathbb{R}^m)  & \rightarrow AC(\mathbb{R}_+,\mathbb{R}^p); \quad
		u & \mapsto &~ Y_{\Sigma,z}[u]
	\end{align*}
	defined by letting $t\mapsto X_{\Sigma,z}[u](t)$ be the solution to the Cauchy problem \eqref{sdyn} with 
  $x_0=z$, and letting $Y_{\Sigma,z}[u](t)=C X_{\Sigma,z}[u](t)$ as in \eqref{out}.
  
  \subsection{Realization of bilinear systems}
  Below, we recall elements of realization theory for bilinear systems \cite{Rugh81}.
  The bilinear system  $\Sigma$ is a \emph{realization} of a function
  \begin{equation} 
  \label{iofunction}
    f: L_{loc}(\mathbb{R}_+,\mathbb{R}^m)  \rightarrow AC(\mathbb{R}_+,\mathbb{R}^p),
  \end{equation}
  if $f=Y_{\Sigma,x_0}$. The system $\Sigma$ is a minimal realization of $f$, if it has the smallest state-space dimension among all bilinear systems which are realizations of $f$. 
  Bilinear system $\Sigma$ is \emph{observable}, if for any two states $x_1,x_2 \in \mathbb{R}^n$, $Y_{\Sigma,x_1}=Y_{\Sigma,x_2}$ implies
  $x_1=x_2$. We say that $\Sigma$ is \emph{span-reachable}, if $\mathrm{Span}\{X_{\Sigma,x_0}[u](t) \mid t \ge 0,u  \in L_{loc}(\mathbb{R}_+,\mathbb{R}^m) \}=\mathbb{R}^n$. 
  It is well-known \cite{Isi:Nonlin} that $\Sigma$ is a minimal realization of $f$, if and only if 
 $\Sigma$ is a realization of $f$, and it is span-reachable and observable. 
 Span-reachability and observability have algebraic characterizations. Specifically, $\Sigma$ is observable, if and only if
$$\bigcap_{w\in Q^*}\ker CA_{w}= \{0\};$$ 
and $\Sigma$ is
 span-reachable if and only if 
 \begin{align*}
  \mathrm{Span}  \{ A_{w}x_0 \mid w \in Q^*\}
 =\mathbb{R}^n.	
 \end{align*}
 In the sequel, we use the notion of generating series for bilinear systems \cite{Isi:Nonlin,Rugh81,GrayPaper}. 
 \begin{definition} 
  A \emph{generating series (over $Q$)} is a function
 $c: Q^{*} \rightarrow \mathbb{R}^{p}$ such that there exist 
 $K,R > 0$ which satisfy 
 \begin{equation}
 	\label{Eq:GrowthCondition}
 	\forall w \in Q^{*}: \|c(w)\| \le KR^{|w|}, 
 \end{equation}
 where  $|w|$ is the length of the sequence $w$ (number of elements from $Q$ in $w$).
 \end{definition}
 
 For each   $u \in L_{loc}(\mathbb{R}_+,\mathbb{R}^m)$, $w \in Q^{*}$ and $t \ge 0$, define the iterated integral
 $V_{w}[u](t)$ as follows: 
 \begin{enumerate}
 	\item $V_{\epsilon}[u](t)=1$,
 	\item for $w=q_1\cdots q_k \in Q^*$, $q_1,\ldots,q_k \in Q$, and $k > 1$ 
 \begin{equation*}
     V_w[u](t)= \int_{0}^{t} u_{q_k}(\tau)V_{q_1\cdots q_{k-1}}[u](\tau) d\tau,
 \end{equation*}
 where $u_0(t)=1$.  
 \end{enumerate}
 For $q \in Q$ we note that $V_{q}[u](t)=\int_0^{t} u_{q}(\tau)d\tau$ since $q=q\epsilon$.
We define the function $F_{c}:L_{loc}(\mathbb{R}_+,\mathbb{R}^m) \rightarrow L_{loc}(\mathbb{R}_+,\mathbb{R}^p)$ generated by a \emph{ generating series} $c$ by
 \begin{equation}
 	\label{Eq:F_functional}
 	F_c[u](t)=\sum_{w \in Q^*} c(w)V_w[u](t). 
 \end{equation} 
 %
 

The map $F_c$ is called a Fliess operator and the right-hand side of \eqref{Eq:F_functional} a Fliess series. The Fliess operator $F_c$ is well-defined, as the growth condition~\eqref{Eq:GrowthCondition} is sufficient for absolute convergence of the series in \eqref{Eq:F_functional}, by Theorem~3.1,~\cite{GrayPaper}.

  Following \cite{Isi:Nonlin}, $f$ has a realization by 
  a bilinear system, only if $f=F_{c_f}$ for some generating series $c_f$. 
  Moreover, a bilinear system $\Sigma$ of the form \eqref{eq:BSSform} is a 
  realization of $f$, if and only if 
  \( \forall w \in Q^{*}: c_f(w)=CA_wx_0. \)


\section{Nice selection for bilinear systems}
\label{sect:nice:def}
At the outset, we define a concept of nice selections. 


\begin{Definition}[Nice selections] \label{def:nice_select}
A subset $\alpha$ of $Q^*$ is called a column nice selection of a bilinear system  $\Sigma$ of the form \eqref{eq:BSSform}, if $\alpha$ has the following property, which we refer to as \emph{prefix closure}: if $w q \in \alpha$ for some $q \in Q$, $w \in Q^*$, then $w \in \alpha$. 
A subset $\beta$ of $Q^{*}$ is called a row nice selection of $\Sigma$ 
if $\beta$ has the following property, which we refer to as  \emph{suffix closure}: if $q w \in \beta$ for some $q \in Q$, $w \in Q^*$, then $w \in \beta$.
\end{Definition}
The phrase ''nice selection`` will be used when it is irrelevant whether a nice row or column selection is used.      

Based on a nice selection, we next introduce the notion of a $\gamma$-partial realizations of an input-output map by $\Sigma$.

\begin{Definition}[$\gamma$-partial realization]
\label{Def:PartialRealisation}
Let $\gamma$ be a nice selection of a bilinear system $\Sigma$ of the form \eqref{eq:BSSform} and let
$f$ be an input-output map of the form \eqref{iofunction}, where it is assumed that $f=F_{c_f}$ for a generating series $c_f$. The bilinear system $\Sigma$ is then called a $\gamma$-partial realization of $f$, if
 $$\forall w \in \gamma : c_f(w)=CA_w x_0.$$
\end{Definition}

Definition \ref{Def:PartialRealisation} states that $\Sigma$ is an $\gamma$-partial realization of $f$, if the generating series of $Y_{\Sigma,x_0}$ and that
of $f$ coincide on the set $\gamma$.

\begin{Remark}[Classical partial realization]
\label{Npart}
Let $\gamma$ be the set of sequences of length at most $N$, $\gamma=\{v \in Q^{*} \mid |v| \le N\}$. In this case we note that
$\gamma$-partial realization reduces to $N$-partial realization as defined in \cite{isi:tac,Isi:Nonlin}.  
In particular, if $N \ge \max\{n,n'\}$, where $n = \dim \Sigma$ and $n'$ is the dimension of a minimal bilinear realization of  $f$ and $\Sigma$ is a $\gamma$-partial realization of 
$f$, then $\Sigma$ is a realization of $f$. 
\end{Remark}

\begin{Remark}[Relationship with moment matching]\label{comp}
As it\\ was mentioned in the introduction, moment matching for 
 bilinear systems was investigated in \cite{Bai2006406,BilinearMomentMatching1,BilinearMomentMatching2,BilinearMomentMatching5}. Below,  we will explain in detail
 the relationship between the existing definition for moment matching and Definition~\ref{Def:PartialRealisation} for $\gamma$-partial realization. 
 Note that in the cited papers, the systems of the form $\dot z = Az+\sum_{i=1}^{m} N_i z u_i+Bu$, $y=Hz$, $z(0)=0$ were studied. By defining $x=(z^T,1)^T$,
 $A_0 = \begin{bmatrix} A & 0 \\ 0 & 0 \end{bmatrix}$, $A_i=\begin{bmatrix} N_i & B_i \\ 0 & 0 \end{bmatrix}$, $i=1,\ldots,m$, $C=\begin{bmatrix} H & 0 \end{bmatrix}$,
 where $B_i$ is the $i$th column of $B$; it is clear that $x$ and $y$ satisfy \eqref{eq:BSSform}. Hence, the system class considered in those paper can be embedded
 into the system class considered in this paper.
 We will consider $m=1$ in order to avoid excessive notation.
 Consider a bilinear system $\Sigma=(p,1,n,\{A_0,A_1\},C,x_0)$ and assume that $\Sigma$ is a realization of $f$. 
 In \cite{BreitenDamm,BennerSiam,FlagGugercin,Wang20121231} the moments 
 $m_{\Sigma}(l_1,\ldots,l_k)$,  $0 < l_1,\ldots,l_k \in \mathbb{N}$ of $\Sigma$
 at certain frequencies $\sigma_1,\ldots,\sigma_k$ were defined as follows: 
  $m_{\Sigma}(l_1,\ldots,l_k)=C(\sigma_1 I - A_0)^{-l_1}A_1(\sigma_2 I - A_0)^{-l_2} \cdots A_1 (\sigma_k I - A_0)^{-l_k}x_0$ for all if $\sigma_1,\ldots, \sigma_k \in \mathbb{C}$ and 
 $m_{\Sigma}(l_1,\ldots,l_k)= c_f(0^{l_1-1}10^{l_2-1} \cdots 10^{l_k-1})=CA_0^{l_1-1}A_1 \cdots  A_1A_0^{l_k-1}x_0$ if $\sigma_1=\cdots=\sigma_k=\infty$. 
 In the cited literature, a system $\bar{\Sigma}=(p,1,n,\{\bar{A}_0,\bar{A}_1\},\bar{C},\bar{x}_0)$
  was said to match the moments of $\Sigma$ for $l_1,\ldots,l_j \in \{0,\ldots,N\}$, $j=1,\ldots,k$, $N \in \mathbb N$ at frequencies
 $\sigma_1,\ldots,\sigma_k$, if $m_{\Sigma}(l_1,\ldots,l_j)=m_{\bar{\Sigma}}(l_1,\ldots,l_j)$ for all  $l_1,\ldots,l_j \in \{1,\ldots,N\}$, $j=1,\ldots,k$.
 Note that in  \cite{Bai2006406,BilinearMomentMatching1,BilinearMomentMatching2,BilinearMomentMatching5} only moments for $\sigma_1=\cdots = \sigma_k=0$ were considered. 
 Moment matching at $\infty$ can be expressed in our framework as follows: $\bar{\Sigma}$ matches the moments of $\Sigma$ $l_1,\ldots,l_j \in \{1,\ldots,N\}$, $j=1,\ldots,k$ at the 
 frequency $\sigma_1=\cdots =\sigma_k=\infty$,  if and only if $\bar{\Sigma}$ is a $\gamma$-partial realization of $Y_{\Sigma,x_0}$, where
 $\gamma=\{ 0^{l_1}1\cdots 0^{l_{j-1}} 10^{l_j}  \mid l_1,\ldots,l_j \in \{0,\ldots,N-1\}, j=1,\ldots,k \}$. Note that $\gamma$ is both prefix and
 suffix closed, i.e., it qualifies both for nice row and nice column selection. 
 For other frequencies, the relationship is less obvious, and it remains a topic of future research.
\end{Remark}

There are two descriptions of the fact that $\Sigma$ is an $\gamma$-partial realization of $f$. The first is that
the input-output map $Y_{\Sigma,x_0}$ of $\Sigma$ is an approximation of $f$ in the sense that \emph{for all} inputs $u$,
the outputs $Y_{\Sigma,x_0}(u)$ and $f(u)$ are close to each other in a suitable metric. The other interpretation is that 
\emph{for some} inputs $u$, $Y_{\Sigma,x_0}(u)$ \emph{ equals } $f(u)$.
Below, we characterize both cases. 


Let $\gamma$ be a nice selection of a bilinear system $\Sigma$ of the form \eqref{eq:BSSform} and for any $T \ge 0$,  define
\begin{multline*}
\mathcal{U}_{\gamma,T}=\{ u \in  L_{loc}(\mathbb{R}_+,\mathbb{R}^m) \mid \forall v \in  Q^{*}, v \notin \gamma, t \in [0,T]:\\ V_{v}[u](t)=0 \}. 
\end{multline*}
For $\alpha \in \mathbb{N}$ and $q \in Q$ let $q^{\alpha}$ denote the sequence $qq\cdots q$ obtained
by repeating $q$ $\alpha$-times. If $\alpha=0$, let $q^{0}$ be the empty word.
Define the set
 \begin{multline*}
  L_{\gamma}=\{ v \in Q^{*} \mid v=q_1\cdots q_k,~ q_1,\ldots,q_k \in Q, k > 0\qquad \\ \qquad\mbox{ such that } 
  \forall v_i \in \{0,q_i\}^{*}, i=1,\ldots,k,  v_1v_2\cdots v_k  \in \gamma \},
 \end{multline*}
and recall that $\{0,q_i\}^{*}$ denotes the set of all sequences of $q_i$ and $0$. In particular, $v \in \{0,q_i\}^*$ if and only if
$v$ contains only the symbols $0$ and $q_i$, or, in other words, $v=0^{l_1}q_i^{l_2}\cdots q_i^{l_{r-1}}0^{l_r}$ for some 
$l_1,\ldots,l_r \in \mathbb{N}$, $r > 0$. We will call  $L_{\gamma}$ the \emph{set of sequences consistent} with $\gamma$.
 
Now for any $q \in Q$, let $e_q$ denotes the $q$th standard basis vector of $\mathbb{R}^m$, if $q \in \{1,\ldots,m\}$, and
$e_0=0$. We will say that $u \in  L_{loc}(\mathbb{R}_+,\mathbb{R}^m)$ is consistent with the nice selection $\gamma$ on an interval $[0,T]$, if 
there exist $q_1,\ldots,q_k \in Q$, reals $0 < t_1 < \cdots  < t_{k-1} < t_k =T$ and scalar valued functions 
$u_i \in L_{loc}([t_i,t_{i+1}],\mathbb{R})$, such that
$u(s)=u_i(s)e_{q_i}$, $s \in [t_{i},t_{i+1})$ for all $i=1,\ldots,k$, $u(T)=e_{q_k}$ and $q_1\cdots q_k \in L_{\gamma}$. That is, if $u$ is consistent with $\gamma$ on $[0,T]$, then $[0,T]$ can be divided into a finite number of intervals, and on each interval, at most one
component of $u$ is not zero.

We are now ready to state the main result relating values of $Y_{\Sigma,x_0}$ and $f$.
\begin{Theorem}[Preserving input-output behavior]
\label{theo1}
Let $\gamma$ be a nice selection, of a bilinear system $\Sigma$, containing all zero sequences; $\{0\}^{*} \subseteq \gamma$.
 \begin{description}
 \item{(A)} If $\Sigma$ is an $\gamma$-partial realization of $f$, then 
       \begin{equation}
       \label{theo1:eq}
        \forall u \in \mathcal{U}_{\gamma,T}, t \in [0,T]: Y_{\Sigma,x_0}(u)(t)=f(u)(t)
      \end{equation}
 \item{(B)}
       Any $u$ consistent with $\gamma$ on $[0,T]$ belongs to $\mathcal{U}_{\gamma,T}$, and $0 \in \mathcal{U}_{\gamma,T}$. 
 \end{description}
\end{Theorem}


The theorem above claims that if $\Sigma$ is a $\gamma$-partial realization of $f$, then $f$ and
the input-output map of $\Sigma$ coincide on the set of inputs which are mapped to the zero map by maps $u\mapsto V_w[u]$ indexed by words not in the nice selection $\gamma$. Moreover, the set of such inputs includes piecewise-constant inputs of switching type. Theorem \ref{theo1} mirrors the results of \cite{MertTAC}.

\begin{pf}
\textbf{Part (A)}
  From $u \in \mathcal{U}_{\beta,T}$ it follows that $V_v[u](t)=0$ for all $v \in Q^{*}, v \notin \gamma$, $t \in [0,T]$, and hence
  $$\sum_{v \notin \gamma, v \in Q^{*}} c(v) V_{v}[u](t)=0=\sum_{v \notin \gamma, v \in Q^{*}} c_f(v) V_{v}[u](t),$$ 
  and \( Y_{\Sigma,x_0}(u)(t)=\sum_{v \in Q^{*}} c(v)V_{v}[u](t)=\sum_{v \in \gamma} c(v) V_{v}[u](t)= \sum_{v \in \gamma} c_f(v) V_{v}[u](t)=\sum_{v \in Q^{*}} c_f(v)V_{v}[u](t)=f(u)(t) \).

 \textbf{Part (B)}
   If $u$ is consistent with $\gamma$, then 
	$$u=\mathbf{u}_1\#_{\tau_1} \mathbf{u}_2 \#_{\tau_2} \cdots \#_{\tau_k} \mathbf{u}_k,$$ 
	where for $i=1,\ldots,k$ 
	$$\tau_1=t_1,\quad\tau_i=t_i-t_{i-1},\quad\mathbf{u}(t)=u_i(t)e_{q_i},$$ 
	and where for any two functions $f_1,f_2$, 
  \[ 
	(f_1 \#_{\tau} f_2)(s)=\left\{\begin{array}{rl}
                             f_1(s) & s \in [0,\tau) \\
                             f_2(s-\tau) & s \in [\tau,+\infty)
                           \end{array}\right.  
\]																							
   Using \cite[eq. (11)]{WangGenSer} repeatedly, it then follows that for any $v \in Q^{*}$ and for any $t \in [t_{j-1},t_j)$, $j=1,\ldots,k$,
   \begin{equation} 
   \label{proof:partB:eq1}
     V_v[u](t)=\sum_{\overset{v_1,\ldots,v_j \in Q^{*},}{ v=v_1v_2\cdots v_j}} V_{v_1}[\mathbf{u}_1](\tau_1) \cdots V_{v_j}[\mathbf{u}_j](\tau_j-t), 
   \end{equation}
   Note that $\mathbf{u}_i(s)=u_i(s)e_{q_i}$ for all $s \in [0,\tau_i]$ and hence $V_{w}[\mathbf{u}_i](s) \ne 0 \implies w \in \{0,q_i\}^{*}$ for all $i=1,\ldots,k$. Hence, from this and \eqref{proof:partB:eq1} it follows that $V_v[u](t) \ne 0$ implies that 
   $v=v_1\cdots v_j$ such that $v_i \in \{0,q_i\}^{*}$, $i=1,..,j$. If $j < k$, we can take $v_{j+1}= \cdots =v_k=\epsilon$, and then 
   $v_i \in \{0,q_i\}^{*}$, $i=1,\ldots,k$,  $v=v_1v_2\cdots v_k$. 
   But by the assumption that $u$ is consistent with $\gamma$ on $[0,T]$, it follows that $q_1\cdots q_k \in L_{\gamma}$, which implies that
   $v=v_1\cdots v_k \in \gamma$. That is, if $V_{v}[u](t) \ne 0$ for some $t \in [0,T]$, then $v \in \gamma$. Hence,
   for any $v \notin \gamma$, $V_v[u](t) = 0$ for all $t \in [0,T]$. The latter implies that $u \in \mathcal{U}_{\gamma,T}$.
\end{pf}

Theorem~\ref{theo1} provides the set of inputs $\mathcal{U}_{\gamma,T}$ such that the $\gamma$-partial realisation and the original input-output map coincides. It is desirable that the set $\mathcal{U}_{\gamma,T}$ is as big as possible; however, the price is a high dimension of the reduced system. Therefore, we use a nice selection to keep the dimension of the reduced system low by pinpointing the attention on specific input trajectories - switching inputs of form  $u(s)=u_i(s)e_{q_i}$ for $s$ in some time interval $ [t_{i},t_{i+1})$.

\begin{Example}
 Let us take $m=2$ and let $\gamma$ to be the set of all sequences containing only $0$'s and $1$'s. 
 The set $\mathcal{U}_{\gamma,T}$ will contain all those input signals whose second component is zero on $[0,T]$, i.e.,
 $u=(u_1,u_2) \in \mathcal{U}_{\gamma,T}$ $\iff$ $\forall t \in [0,T]: u_2(t)=0$. 

 Let us now take   $\gamma=\{v_1v_2  \mid v_1 \in \{0,1\}^{*}, v_2 \in \{0,2\}^{*} \}$. 
Then 
 for any $T,t_1 > 0$, the input 
 \( u(t)=\left\{\begin{array}{rl} (1,0)^T & t < t_1 \\
                                  (0,1)^T & t_1 \le t
                \end{array}\right.
  \)                                          
  will belong to $\mathcal{U}_{\gamma,T}$,
   but 
  \(
u(t)=\left\{\begin{array}{rl} (0,1)^T & t < t_1 \\
                                  (1,0)^T & t_1 \le t
                \end{array}\right.
  \)
 will not belong to $\mathcal{U}_{\gamma,T}$, since $V_{21}[u](t)= \int_0^{T} u_1(\tau)\int_0^{\tau} u_2(\tau_1)d\tau_1d\tau=\int_{t_1}^{T} \int_0^{t_1} d\tau_1d\tau=t_1(T-t_1) \ne 0$, but the sequence $21 \notin \gamma$. 
\end{Example}


 Let us now investigate the behavior of $\gamma$-partial realization of $f$ for all bounded inputs.   
 We formulate the following simple but useful observation, where for any $\gamma\subseteq Q^*$ containing the empty sequence $\epsilon$, we introduce the notation 
\begin{align}
\label{Eq:MaximalLength}
	N_\gamma=\max\{ N\in\mathbb{N} ~|~ \{ v \in Q^{*} \mid |v| \le N\} \subseteq \gamma\},
\end{align}
which in ''worst`` case is zero.
 \begin{Lemma}
 \label{Lemma:BondsOnC}
 	Let $c_f$ be the generating series of $f$ and  $c$ be the generating series of $Y_{\Sigma,x_0}$. 
 	Then there exist real numbers $K_{\gamma}, M_{\gamma} > 0$ such that 
\begin{equation}
   \label{lemma1:eq}
    \forall w \in Q^{*} \setminus \gamma: \| c_f(w)-c(w) \| \le K_{\gamma} M_{\gamma}^{|w|-N_\gamma}. 
   \end{equation}
 \end{Lemma}
 
 \begin{pf}
 	Recall that $c(w) = C A_w x_0$. Since $c_f$ is a generating series, there exists $K_f > 0, M_f > 0$ such that for all $w \in Q^{*}: \| c_f(w) \|  < K_f M_f^{|w|}$.
 	Leaning on the two observations above, we define 
 	 \(	M_A =\max\{\|A_0\|,\ldots, \|A_m\|\} \), 
 		\( K'_{\gamma} = \max\{\|C\| M_A^{N_\gamma} \|x_0\|, K_f M_f^{N_\gamma}\} \),
		\( M_{\gamma} =\max\{M_f,M_A\} \), and 
		\( K_\gamma=2K'_\gamma \).
 \end{pf}


Lemma~\ref{Lemma:BondsOnC} is an ingredient in the proof of Theorem~\ref{theo2}, which provide an error bound between the input output map $f$ and a corresponding $\gamma$-partial realization, on the input trajectories not belonging to $\mathcal{U}_{\gamma,T}$.

 \begin{Theorem}
 \label{theo2}
  Let $\gamma$ be a nice selection and suppose that  $u \in  L_{loc}(\mathbb{R}_+,\mathbb{R}^m)$ and 
  \[ \sup_{t \in [0,T]} \|u(t)\| < R, \]
   and $\Sigma$ is a 
  $\gamma$-partial realization of $f$. Then  there exist 
	reals $K_{\gamma}, M_{\gamma} > 0$ such that
  \begin{equation}
\label{error:bound}
  \begin{split}
  & \forall t \in [0,T]: \| Y_{\Sigma,x_0}(u)(t) - f(u)(t)\|  \\ 
  & \leq
 K_{\gamma}((m+1)^2 \max\{R,t\})^{N_{\gamma}} e^{M_{\gamma}(m+1)^2\max\{R,t\}}
 \end{split}
\end{equation}
 \end{Theorem}

%
%

The theorem 
above complies with the intuition that for sufficiently small time and small inputs, i.e., for $(m+1)^2 \max\{R,t\} < 1$, the larger $N_{\gamma}$ the smaller is the difference between $f$ and the input-output map of $\Sigma$. 

\begin{pf}
By Lemma~2.1 in \cite{GrayPaper}, if $K=\max\{R,t\}$ then for any $w = q_1,\ldots,q_k \in Q^*$, 
\begin{equation} 
\label{bound}
   \|V_{w}[u](t)\| \le K^{k} \frac{1}{r_0! \cdots r_m!} 
\end{equation}
where $r_i \equiv r_i(w)$ is the number of occurrences of the integer $i \in Q$ in the sequence $w$. We will use a  notation $r(w) = (r_0(w), \hdots, r_m(w))$.
To prove~\eqref{error:bound}, we use \eqref{lemma1:eq} and the following observation. By the binomial expansion, we have 
\(
(m+1)^k = ((1 + \hdots +1) +1)^k = \sum_{|\bar{r}|=k} \frac{k!}{\bar{r}_0! \hdots \bar{r}_m!}.
\)
where we have used the notation $|\bar{r}| \equiv \bar{r}_0+\hdots+\bar{r}_m$.

As a consequence,
\[
  \begin{split}
 & \| Y_{\Sigma,x_0}(u)(t) - f(u)(t)\| \le  \sum_{w \in Q^{*} \setminus \gamma} \| c_f(w)-c(w)\| |V_w[u](t)|  \\
 & \le  \sum_{w \in Q^{*} \setminus \gamma} K_{\gamma} M_{\gamma}^{|w|-N_{\gamma}} \frac{K^{|w|}}{r_0(w)! \cdots r_m(w) !}\\
    \end{split}
\]
for $N_{\gamma}$ defined in \eqref{Eq:MaximalLength}. We continue 
  \[
  \begin{split}
 &\sum_{w \in Q^{*} \setminus \gamma} K_{\gamma} (M_{\gamma})^{|w|-N{\gamma}} \frac{K^{|w|}}{r_0(w)! \cdots r_m(w) !} \\
 & \le \sum_{w \in Q^{*}, |w| \ge N_{\gamma}} K_{\gamma} (M_{\gamma})^{|w|-N_\gamma} \frac{K^{|w|}}{|w|!} \frac{|w|!}{r_0(w)! \cdots r_m(w) !} \\
 &= \sum_{k=N_{\gamma}}^{\infty} K_{\gamma} (M_{\gamma})^{k-N_{\gamma}} \frac{K^k}{k!} 
  \sum_{\tiny\begin{matrix}
  	|\bar{r}|=k \\ \bar{r} = (\bar{r}_0, \hdots, \bar{r}_m)
  \end{matrix}}
  \sum_{
  \tiny\begin{matrix}
  	w \in Q^{*} \\
  	r(w)=\bar{r}
  \end{matrix}
   } \frac{k!}{\bar{r}_0! \cdots \bar{r}_m!} \\
  & \le  \sum_{k=N_{\gamma}}^{\infty} K_{\gamma} (M_{\gamma})^{k-N_{\gamma}} \frac{K^k}{k!}  \left(\sum_{\bar{r}_0+\cdots +\bar{r}_m=k} \frac{k!}{\bar{r}_0! \cdots \bar{r}_m!} \right)^2 \\ 
  & = \sum_{k=N_{\gamma}}^{\infty} K_{\gamma} (M_{\gamma})^{k-N_{\gamma}} \frac{K^k}{k!} (m+1)^{2k} \\
 &= K_{\gamma} (K (m+1)^{2})^{N_{\gamma}}   \sum_{k=0}^{\infty} \frac{1}{(k+N_{\gamma})!} (M_{\gamma})^{k} (K (m+1)^{2})^{k}  \\
 & \le  K_{\gamma} (K (m+1)^{2})^{N_{\gamma}} \sum_{k=0}^{\infty} \frac{1}{k!} M_{\gamma}^{k} (K (m+1)^{2})^{k} = \\
 & = K_{\gamma}((m+1)^2K)^{N_{\gamma}} e^{M_{\gamma}(m+1)^2K}
    \end{split}
 \]
We have used the observation that  
%
for $\bar{r}=(\bar{r}_0, \hdots, \bar{r}_m),$ with $|\bar{r}|=k$, $|\{w \in Q^{*}|~ r(w)=\bar{r}\}| = \frac{k!}{\bar{r}_0! \cdots \bar{r}_m!}$.

\end{pf}

In summary, Theorem~\ref{theo1} characterises the set $\mathcal{U}_{\gamma,T}$ of input trajectories $U$ for which input output map of the original and reduced system are identical; whereas, Theorem~\ref{theo2} provides the error bounds for the inputs which do not belong to $\mathcal{U}_{\gamma,T}$.

\section{Model reduction by nice selection}
\label{Sec:Algorithms}
In this section, we present procedures for computing $\alpha$-partial and $\beta$-partial realizations of an input-output map $f$ which is realizable by a bilinear system.

\begin{Definition}
Let  $\Sigma$ be a bilinear systems of the form \eqref{eq:BSSform}. Let $\alpha$ be a nice row selection and $\beta$ be a nice column selection related to $\Sigma$. Then the subspaces
\begin{equation*}
\begin{aligned}
& \mathscr{O}_{\alpha}(\Sigma)= \bigcap_{w \in \alpha} \ker C A_w  \hbox{ and }\\
& \mathscr{R}_{\beta}(\Sigma) = \SPAN \{ A_wx_0 \mid w \in \beta\}
\end{aligned}
\end{equation*}
will be called \emph{$\alpha$-unobservability and $\beta$-reachability spaces} of $\Sigma$ respectively.
\end{Definition}
The spaces $\mathscr{O}_\alpha(\Sigma)$ and $\mathscr{R}_\beta(\Sigma)$ will be denoted by $\mathscr{O}_\alpha$ and $\mathscr{R}_\beta$,  if $\Sigma$ is clear from the context. 
\begin{Theorem} \label{theo:krylov1}
	Let $\Sigma=(p,m,n, \{A_i\}_{i\in Q}, C ,x_0)$  be a bilinear system, and $\beta$ be a nice column selection.  Let $V \in \mathbb{R}^{n \times r}$ be a full column rank matrix such that
	\[
	\mathscr{R}_{\beta} = \IM (V), 
	\]
	and let $V^{-1}$ be any left inverse of $V$. Define 
	\[
	\begin{split}
	\forall q \in Q: \bar{A}_q=V^{-1}A_qV, \quad \bar{C} =C V, \quad \bar{x}_0=V^{-1}x_0.
	\end{split}
	\]
	Then $\bar{\Sigma}=(p,m,r, \{\bar{A}_q\}_{q\in Q}, \bar{C} ,\bar{x}_0)$
	is a $\beta$-partial realization of $f=Y_{\Sigma,x_0}$. Furthermore, the $\beta$-reachability spaces of $\bar{\Sigma}$ and $\Sigma$ are equal, $\mathscr{R}_\beta(\bar{\Sigma}) = \mathscr{R}_\beta(\Sigma)$.  
\end{Theorem}
The proof of Theorem~\ref{theo:krylov1} follows the idea of the proofs of Theorems 3 and 6 in \cite{MertTAC} and is omitted.  
By duality, we can formulate moment matching by nice row selections, as in Theorems \ref{theo:krylov2}.
\begin{Theorem} \label{theo:krylov2}
	Let $\Sigma=(p,m,n, \{A_i\}_{i\in Q}, C ,x_0)$  be a bilinear system and let $\alpha$ be a nice row selection. Let $W \in \mathbb{R}^{r \times n}$ be a full row rank matrix such that
	\[
	\mathscr{O}_\alpha = \ker (W), 
	\]
	and let $W^{-1}$ be any right inverse of $W$.  
	\[
	\begin{split}
	\forall q \in Q: \bar{A}_q=WA_qW^{-1} \mbox{, } \bar{C} =C W^{-1} \mbox{, and } \bar{x}_0=Wx_0.
	\end{split}
	\]
    Then $\Sigma=(p,m,r, \{\bar{A}_q\}_{q\in Q}, \bar{C} ,\bar{x}_0)$ 
    is an $\alpha$-partial realization of $f=Y_{\Sigma,x_0}$. Furthermore, the $\alpha$-unobservability spaces of $\bar{\Sigma}$ and $\Sigma$ are equal, $\mathscr{O}_\alpha(\bar{\Sigma}) = \mathscr{O}_\beta(\Sigma)$. 
\end{Theorem}
\begin{Remark}[Choosing the model order]
\label{choice:order}
 Note that using nice selections allows us to choose the order of
 the reduced system. Indeed, assume that $\gamma$ is a nice row or column selection, and assume that $\gamma$ has $\bar{n}$ elements. It is then easy to see that $\dim \mathscr{R}_{\gamma} \le \bar{n}$. If $p=1$, i.e. there is one output, then
$n-\dim \mathscr{O}_{\gamma} \le  \bar{n}$.
Moreover, in this case, 
generically, $\dim \mathscr{R}_{\gamma}=\bar{n}$ and $\dim \mathscr{O}_{\gamma}=n-\bar{n}$. This can be shown in a manner similar to
the the discussion after \cite[Theorem 7]{MertTAC}.
Hence, Theorem \ref{theo:krylov1} -- \ref{theo:krylov2} yield
a reduced order model of order at most $\bar{n}$. Similarly to the proof of \cite[Theorem 7]{MertTAC}, it can be shown that
if  $\Sigma$ is a minimal system, then for any $r=1,\ldots,n$, there exist a 
nice row selection $\alpha$ and a nice column selection $\beta$ such that $\alpha$ and $\beta$ both have $r$ elements and 
$\dim \mathscr{R}_{\beta}=r$, $\dim \mathscr{O}_{\alpha}=n-r$, and Theorem \ref{theo:krylov1} -- \ref{theo:krylov2} yield
a reduced order model of order $r$.  In Section \ref{num}, we illustrate this in a numerical example. 
The discussion above can be extended to the case with multiple outputs, if the definition of a nice row selection is modified to include the choice of the
output channel. This can be done along the lines of \cite[Definition 4]{MertTAC}. 
\end{Remark}
If the nice selections at hand are finite and of small size, then  computing matrix representations of the spaces $\mathscr{R}_{\beta}$ and $\mathscr{O}_{\alpha}$ is trivial.
However,  nice row and column selections need not to be finite, or if they are finite, their cardinality can be large. For example, the nice selections which satisfy 
the conditions of Theorem \ref{theo1} are always infinite, and the nice selection corresponding to classical $N$-partial realization (see Remark \ref{Npart}) is
finite but its cardinality is exponential in $N$. Hence, for these cases, the question arises how to compute matrices $V$ and $W$ used in Theorem \ref{theo:krylov1} -- \ref{theo:krylov2}. 

We will start by presenting a special case, when $\alpha=\{ v \in Q^{*} \mid |v| \le N\}$ and $\beta = \{ v \in Q^{*} \mid |v| \le N\}$. 
Denote $\mathscr{R}_{\beta}$ by $\mathscr{R}_N$ and 
$\mathscr{O}_{\alpha}$ by $\mathscr{O}_N$.
Denote by $\mathbf{orth}(M)$ the orthogonal matrix $V$ such that $V$ is full column rank, $\IM(V)=\IM (M)$ and $V^TV=I$. 
Consider the algorithms Algorithm \ref{alg1} -- \ref{alg2}.

\begin{algorithm}
	\caption{
	Calculate  a matrix representation of $\mathscr{R}_N$,
		\newline
		\textbf{Inputs}: $(\{A\}_{q \in Q},x_0)$ and $N$
		\newline
		\textbf{Outputs:} $V  \in \mathbb{R}^{n \times r}$, $\Rank (V)=r$,
		$\IM (V) = \mathscr{R}_N$.
	}
	\label{alg1}
\begin{algorithmic}[1]
		\STATE $V_0:=U_0$; $U_0:=x_0$
		\FOR{$k=1\ldots N$} 
		\STATE
		$V:=\mathbf{orth}(\begin{bmatrix} U_0, & A_0V, & A_1V, & \ldots, & A_{m}V \end{bmatrix})$
		\ENDFOR
		\RETURN $V$.
	\end{algorithmic}
\end{algorithm}
\begin{algorithm}
	\caption{
		Calculate a matrix representation of $\mathscr{O}_N$
		\newline
		\textbf{Inputs}: $(C,\{A_q\}_{q \in Q})$ and $N$
		\newline
		\textbf{Output:} $W \in \mathbb{R}^{r \times n}$, $\Rank (W) = r$ and $\ker (W)=\mathscr{O}_N$.
	\label{alg2}
	}
	\begin{algorithmic}[1]
        	\STATE $V:=U_0$, $U_0:=C^T$.
		\FOR{$k=1\ldots N$} 
		\STATE
		$V:=\mathbf{orth}(\begin{bmatrix} U_0, & A_0^TV, & A_1^TV, & \ldots, & A_{m}^TV \end{bmatrix})$
		\ENDFOR
		\RETURN $W=V^{\mathrm{T}}$.
	\end{algorithmic}
\end{algorithm} 
\begin{Lemma}
  With the notation above, Algorithm \ref{alg1} returns a matrix $V$ such that $\IM (V) = \mathscr{R}_N$ and $V^TV=I_r$ and 
  Algorithm \ref{alg2} returns a matrix $W$ such that $\ker W=\mathscr{O}_N$.
\end{Lemma}
Notice that the computational complexity of Algorithm \ref{alg1} and Algorithm \ref{alg2} is polynomial in $N$ and $n$, even though the spaces of $\mathscr{R}_N$ (resp. $ \mathscr{O}_N$) are generated by images (resp. kernels) of exponentially many matrices.
 In fact, Algorithm \ref{alg1} -- \ref{alg2} stop after $N$ iterations.

For more general nice row and column selections, we need more sophisticated algorithms. To this end,  we recall the concept of non-deterministic finite state automaton and its language.
\begin{Definition} \label{def:NDFA}1
	A non-deterministic finite state automaton (NDFA) is a tuple $\mathcal{A}=(S,Q,\{\rightarrow_q\}_{q \in Q} ,F, s_0)$ such that
	\begin{enumerate}
		\item $S$ is the finite state set,  
                \item $F \subseteq S$ is the set of accepting (final) states,
		\item $\rightarrow_q \subseteq S \times S$ is the state transition relation labelled by $q$, and
		\item $s_0 \in S$ is the initial state.
	\end{enumerate}
	For every $w \in Q^{*}$, define $\rightarrow_w$ inductively as follows: $\rightarrow_{\epsilon}=\{ (s,s) \mid s \in S \}$ and $\rightarrow_{wq} = \{ (s_1,s_2) \in S \times S \mid \exists s_3 \in S: (s_1,s_3) \in \rightarrow_{w} \mbox{ and } (s_3,s_2) \in \rightarrow_q \}$ for all $q \in Q$. We denote the fact $(s_1,s_2) \in \rightarrow_w$ by $s_1 \rightarrow_w s_2$. Define the language $L(\mathcal{A})$ accepted by $\mathcal{A}$ as 
	\[
	L(\mathcal{A})=\{ w \in Q^{*} \mid \exists s \in F: s_0 \rightarrow_w s \}.
	\]
	We say that $\mathcal{A}$ is \emph{co-reachable}, if from any state a final state can be reached, i.e., for any $s \in S$, there exists $w \in Q^*$ and $s_f \in F$ such that $s \rightarrow_w s_f$. It is well-known that if $\mathcal{A}$ accepts $L$, then we can always compute an NDFA $\mathcal{A}_{co-r}$ from $\mathcal{A}$ such that $\mathcal{A}_{co-r}$ accepts $L$ and it is co-reachable.
\end{Definition}
Consider a nice row selection $\alpha$ and a nice column selection $\beta$, and assume that $\alpha$ and $\beta$ are regular languages.
The matrix representations of the spaces $\mathscr{R}_{\beta}$ and $\mathscr{O}_{\alpha}$ can be computed by Algorithms \ref{alg4} -- \ref{alg5}.
\begin{algorithm}
	\caption{
		Calculate  a matrix representation of $\mathscr{R}_{\beta}$, 
		\newline
		\textbf{Inputs}: $(\{A_q\}_{q \in Q},x_0)$ and $\hat{\mathcal{A}}=(S,\{\rightarrow_q \}_{q \in Q},F,s_0)$ such that $L(\hat{\mathcal{A}})=\beta$, $F=\{s_{f_1},\dots s_{f_k}\}$, $k \geq 1$ and $\hat{\mathcal{A}}$  is co-reachable.
		\newline
		\textbf{Outputs:} $\hat{V}  \in \mathbb{R}^{n \times \hat{r}}$, $\Rank(\hat{V})=\hat{r}$,
		$\IM (\hat{V}) = \mathscr{R}_{\beta}$. 
	}
	\label{alg4}
	\begin{algorithmic}[1]
		\STATE $\forall s \in S \backslash \{s_0\}: V_s:=0$, $V_{s_0}:=\mathbf{orth}(x_0)$. 
		\label{alg4.0}
                \REPEAT
		\label{alg4.1}
		\STATE $\forall s \in S: V_s^{old} := V_s$
		\FOR{$s \in S$}
		\STATE  $M_s:=V_s$
		\FOR{$q \in Q, s^{'} \in S: s^{'} \rightarrow_q s$}
		\STATE $M_s:=\begin{bmatrix} M_s, & A_qV^{old}_{s^{'}} \end{bmatrix}$
		\ENDFOR
		\STATE $V_s := \mathbf{orth}(M_s)$
		\ENDFOR
		\UNTIL{$\forall s \in S: \Rank (V_s) = \Rank (V^{old}_s$)}
		\RETURN $\hat{V}:=\mathbf{orth} \left( \begin{bmatrix} V_{s_{f_1}} & \cdots & V_{s_{f_k}} \end{bmatrix} \right)$.
	\end{algorithmic}
\end{algorithm}
\begin{algorithm}
	\caption{
		Calculate a matrix representation of $\mathscr{O}_{\alpha}$, 
		\newline
		\textbf{Inputs}: $(\{A_q\}_{q \in Q},C)$ and $\hat{\mathcal{A}}=(S,\{\rightarrow_q \}_{q \in Q},F,s_0)$ such that $L(\hat{\mathcal{A}})=\alpha$, $F=\{s_{f_1},\cdots s_{f_k}\}$, $k \geq 1$ and $\hat{\mathcal{A}}$  is co-reachable.
		\newline
		\textbf{Outputs:} $\hat{W}  \in \mathbb{R}^{\hat{r} \times n}$, $\Rank(\hat{W})=\hat{r}$,
		$\ker(\hat{W}) = \mathscr{O}_{\alpha}$.
	}
	\label{alg5}
	\begin{algorithmic}[1]
		\STATE $\forall s \in S \backslash F: W_s:=0$, $\forall s \in F: W_{s}^{\mathrm{T}}:=\mathbf{orth}((H)^{\mathrm{T}})$. 
		\label{alg2.0}
                \REPEAT
		\label{alg2.1}
		\STATE $\forall s \in S: W_s^{old} := W_s$
		\FOR{$s \in S$}
		\STATE  $M_s:=W_s$
		\FOR{$q \in Q, s^{'} \in S: s \rightarrow_q s^{'}$}
		\STATE $M_s:=\begin{bmatrix} M_s^T, &  (W^{old}_{s^{'}}A_q)^T \end{bmatrix}^T$
		\ENDFOR
		\STATE $W_s^{\mathrm{T}} := \mathbf{orth}(M_s^{\mathrm{T}})$
		\ENDFOR
		\UNTIL{$\forall s \in S: \Rank (W_s) = \Rank (W^{old}_s)$}
		\RETURN $\hat{W}:=W_{s_0}$.
	\end{algorithmic}
\end{algorithm}
\begin{Lemma}
\label{lem:correctness} 
 Suppose that $\alpha$ and $\beta$ are regular languages, i.e., there is a co-reachable NDFA $\hat{\mathcal{A}}$ which accepts $\beta$ respectively $\alpha$. Algorithm \ref{alg4} returns a full column rank matrix $\hat{V}$ such that $\IM (\hat{V})= \mathscr{R}_{\beta}$, and Algorithm \ref{alg5} returns a full row rank matrix $\hat{W}$ such that $\ker(\hat{W}) = \mathscr{O}_{\alpha}$.  
\end{Lemma}

Algorithms \ref{alg1} -- \ref{alg2}  are particular instances of Algorithm \ref{alg4} -- \ref{alg5} when applied to the nice selections 
$\alpha=\beta=\{v \in Q^{*} \mid |v| \le N\}$ with the choice of accepting automaton $\bar{\mathcal{A}}=(S,Q,\{\rightarrow_q\}_{q \in Q}, S, 0\}$, where $S=\{0,\ldots,N\}$,
 $\rightarrow_{q}=\{ (i,i+1) \mid i=0,\ldots,N-1\}$ for all $q \in Q$.

 As a final remark, Algorithms \ref{alg4}--\ref{alg5} are essentially the same as  \cite[Algorithm 1 -- 2]{MertTAC} and Lemma \ref{lem:correctness} is a simple consequence of \cite[Lemma 1]{MertTAC}. 
In particular, as it was explained in \cite{MertTAC}, 
the computational complexity of Algorithms \ref{alg4}--\ref{alg5} is polynomial in $n$, and the
algorithms stop after at most $n|S|$ iterations, where $|S|$ is the number of states of the NDFA
$\hat{\mathcal{A}}$. However, $n|S|$ is the upper bound on the number of iterations, for particular
examples the actual number of iterations can be much smaller. 

\section{Example}
\label{num}
 To illustrate the methods developed in this paper, we consider a number of numerical examples of systems.
 
 We start with a small dimensional example with the aim of illustrating Theorem 1 and Lemma 1. Consider the bilinear system $\Sigma$ of the form \eqref{eq:BSSform}, where
\[
  \begin{split}
  & A_0=-\begin{bmatrix} 0 & 0 & 0 & 0 \\
             0 & 0 & 0 & 0 \\
             0 & 0 & 1 & 0 \\
             0 & 0 & 0 & 0
        \end{bmatrix},~
   A_1= \begin{bmatrix} 0 & 0 & 0 & 0 \\
           0 & 0 & 0 & 0 \\
           1 & 0 & 0 & 0 \\
           0 & 0 & 0 & 0
        \end{bmatrix},~
   A_2 =\begin{bmatrix} 0 & 0 & 0 & 10 \\
            0 & 0 & 0 & 0 \\
            0 & 0 & 0 & 0 \\
            0 & 0 & 0 & 0
       \end{bmatrix}, \\ 
  & A_3 =\begin{bmatrix}  0 & 1 & 0 & 0 \\
             -3 & -0.1 &  0 & 0 \\
             0 & 0 & 2 & 0 \\
             0 & 0 &  0 & -1
       \end{bmatrix},~~ 
  C=\begin{bmatrix} 1 \\ 0  \\ 1  \\ 0 \end{bmatrix}^T,~~ 
 x_0=\begin{bmatrix} 0 \\ 0 \\ 0 \\ 1 \end{bmatrix}
\end{split}
\]
 with $m=3$, $n=4$. Consider the following nice column selection
 $\gamma$ accepted by an NDFA of the form $\mathcal{A}=(S,Q,\{ \rightarrow_q \}_{q \in Q},F,s_0)$, where $S=\{1,2,3\}=F$,  $s_0=1$ and
 $\rightarrow_{0} =\{ (i,j) \mid i \le j \} \cup \{ (3,1) \}$ and
 for all $q=1,\ldots,m$, $\rightarrow_{q} =\{ (i,q) \mid i \le q \}$.
  In other words, $\gamma = L(\mathcal{A})$.  More explicitly,
 $\gamma$ can be described as follows. Let $L$ be the set of all
  the sequences $w \in Q^{*}$ such that $w=v_1v_2v_3$, 
  $v_i \in \{0,i\}^{*}$, $i=1,2,3$.
 Then $\gamma$ is the set of words of the form
 $w_1w_2\cdots w_k$, $k \ge 1$, such that 
 $w_1 \in L$ and for all $i=2,\ldots,k$, $w_i=0w_i^{'}$ where 
 $w_i^{'} \in L$. To illustrate, the word $1^*2^*3^*0^*2^*3^* \in \gamma$, whereas $2^{l_1} 1^{l_2} 3^*$ for $l_1, l_2 \in \mathbb{N} \setminus \{0\}$   does not belong to $\gamma$. 
 
\begin{figure}[hbpt]
	\centering
  \includegraphics[scale=0.45]{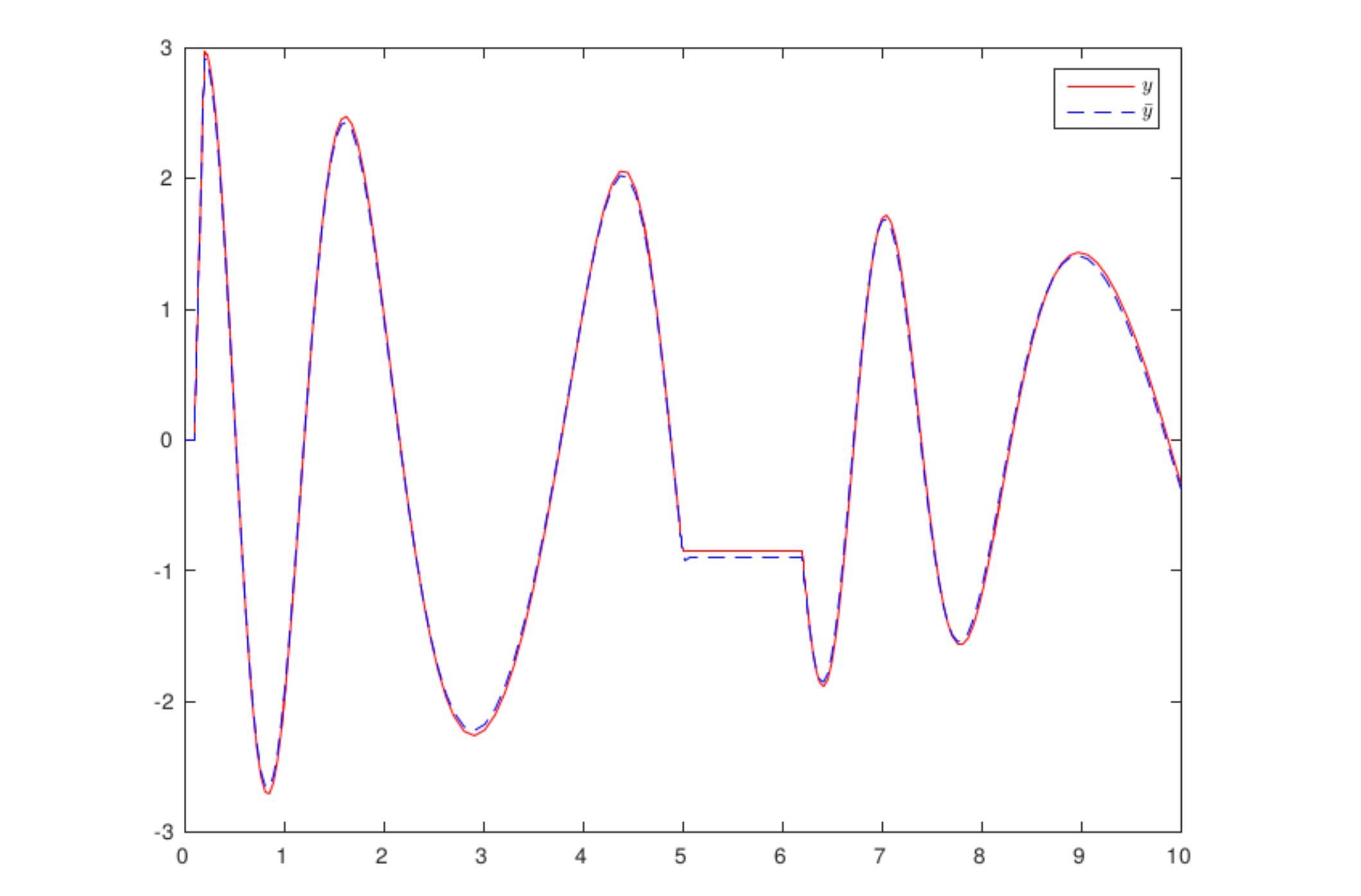}
	\caption{Comparisons of responses of $\Sigma$ (red) and $\bar{\Sigma}$ (green) for $u$ defined in \eqref{input1}. \label{fig1}}
  \end{figure}

 Algorithm \ref{alg4} and Theorem \ref{theo:krylov1} yield
 a reduced order model $\bar{\Sigma}=(\bar{C},\{\bar{A}\}_{q \in Q},\bar{x}_0)$ 
 of size $r=3$, where 
 \[ 
   \begin{split}
    & \bar{A}_0=\bar{A}_1=\mathbf{O}_{3 \times 3}, \quad  \bar{A}_1=\begin{bmatrix} 0 &  10 &       0 \\
         0     &    0  &        0 \\
         0     &    0  &      0
        \end{bmatrix}, ~ 
       \bar{A}_2 = \begin{bmatrix}  0   &     0 &    1  \\
         0 &   -1 &         0 \\
       -3  &   0  &  -0.1
      \end{bmatrix} \\
   & \bar{C}=\begin{bmatrix} 1 & 0 & 0 \end{bmatrix}, \quad 
     \bar{x}_0=\begin{bmatrix} 0 & 1 & 0 \end{bmatrix}^T 
   \end{split}
 \]
 and $\mathbf{O}_{3 \times 3}$ denotes the $3 \times 3$ zero matrix.
 By Theorem \ref{theo:krylov1}, $\bar{\Sigma}$ is a $\gamma$-partial
 realization of the input-output map $f=Y_{\Sigma,x_0}$; and by Theorem \ref{theo1},
  for any $u \in \mathcal{U}_{\gamma,T}$, 
 the outputs of $\bar{\Sigma}$ and $\Sigma$ are equal on $[0,T]$. 
 For example, consider the input 
 \begin{equation}
 \label{input1}
    u(t) = \left\{\begin{array}{rl}
             (\cos(t\pi)+2,0,0)^T & t \in [0,0.1) \\
             (0,\cos(t\pi)+2,0)^T & t \in [0.1,0.2) \\
             (0,0,\cos(t\pi)+2)^T & t \in [0.2,5) \\
             (0,0,0)^T           & t \in [5,6) \\
             (0,\cos(t\pi)+2,0)^T & t \in [6.1,6.2) \\
             (0,0,\cos(t\pi)+2)^T & t \in [6.1,+\infty), \\
            \end{array}\right.
 \end{equation}
  then by Theorem \ref{theo1}, $u \in \mathcal{U}_{\gamma,T}$ for all $T \ge 6.1$ and hence
  $y=Y_{\Sigma,x_0}(u)$ and $\bar{y}=Y_{\bar{\Sigma},\bar{x}_0}(u)$ coincide on $[0,T]$. The corresponding response
  is shown on Figure \ref{fig1} for $T=10$. Although according to Theorem \ref{theo1}, these two responses
should be equal, however, on Figure \ref{fig1} one sees a slight difference, which is due to numerical error.   

\begin{figure}[hbpt]
	\centering
  \includegraphics[scale=.6]{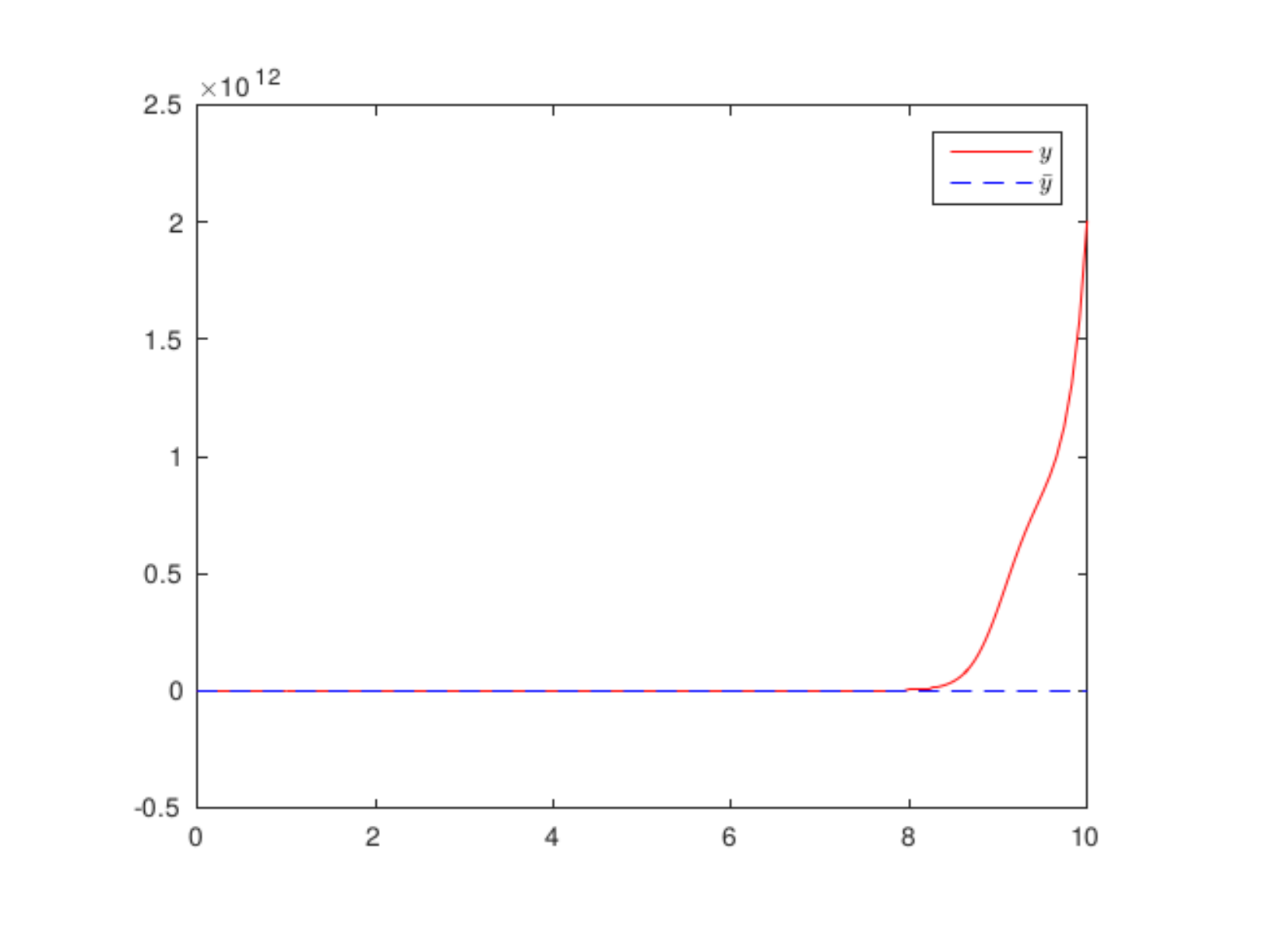}
	\caption{Comparisons of responses of $\Sigma$ (red) and $\bar{\Sigma}$ (blue) for $u$ defined in \eqref{input2}. \label{fig2} }
  \end{figure}

	However, for $u$ given by   
  \begin{equation}
 \label{input2}
   u(t) = \left\{\begin{array}{rl}
             (0,\sin(t\pi)+2,0)^T & t \in [0,0.5) \\
             (\sin(t\pi)+2,0,0)^T & t \in [0.5,1) \\
             (0,0,\sin(t\pi)+2)^T & t \in [1,+\infty) \\
            \end{array}\right.
 \end{equation} 
$u \notin \mathcal{U}_{\gamma,T}$, $T=10$. The corresponding outputs $y=Y_{\Sigma,x_0}(u)$ and $\bar{y}=Y_{\bar{\Sigma},\bar{x}_0}(u)$ need not coincide on $[0,T]$. The corresponding response
  are shown in Figure \ref{fig2} for $T=10$. In fact, $\Sigma$ is exponentially unstable for this $u$, while the state and output
  trajectory of $\bar{\Sigma}$ remain bounded.    
  
  Note that if we choose the nice selection $\alpha=\{ \epsilon, 3, 34\}$, and we apply Theorem \ref{theo:krylov1} to $\Sigma$ above, then $\mathscr{R}_{\alpha}$ will be an 
  invertible $3 \times 3$ matrix and in this case Theorem \ref{theo:krylov1} yields a reduced order bilinear system which is related to $\bar{\Sigma}$ described above by a 
  linear isomorphism. If we apply Theorem \ref{theo:krylov1} to $\alpha=\{\epsilon,3\}$, then we obtain a reduced order model of order $2$. This illustrates that
  finite  nice selections can be used to choose the order of the reduced model, as explained in Remark \ref{choice:order}. 
  The code for this example can be found in the supplementary material, in the file \texttt{BilinearModelExample2.m}.

\begin{figure}[hbpt]
	\centering
  \includegraphics[scale=0.45]{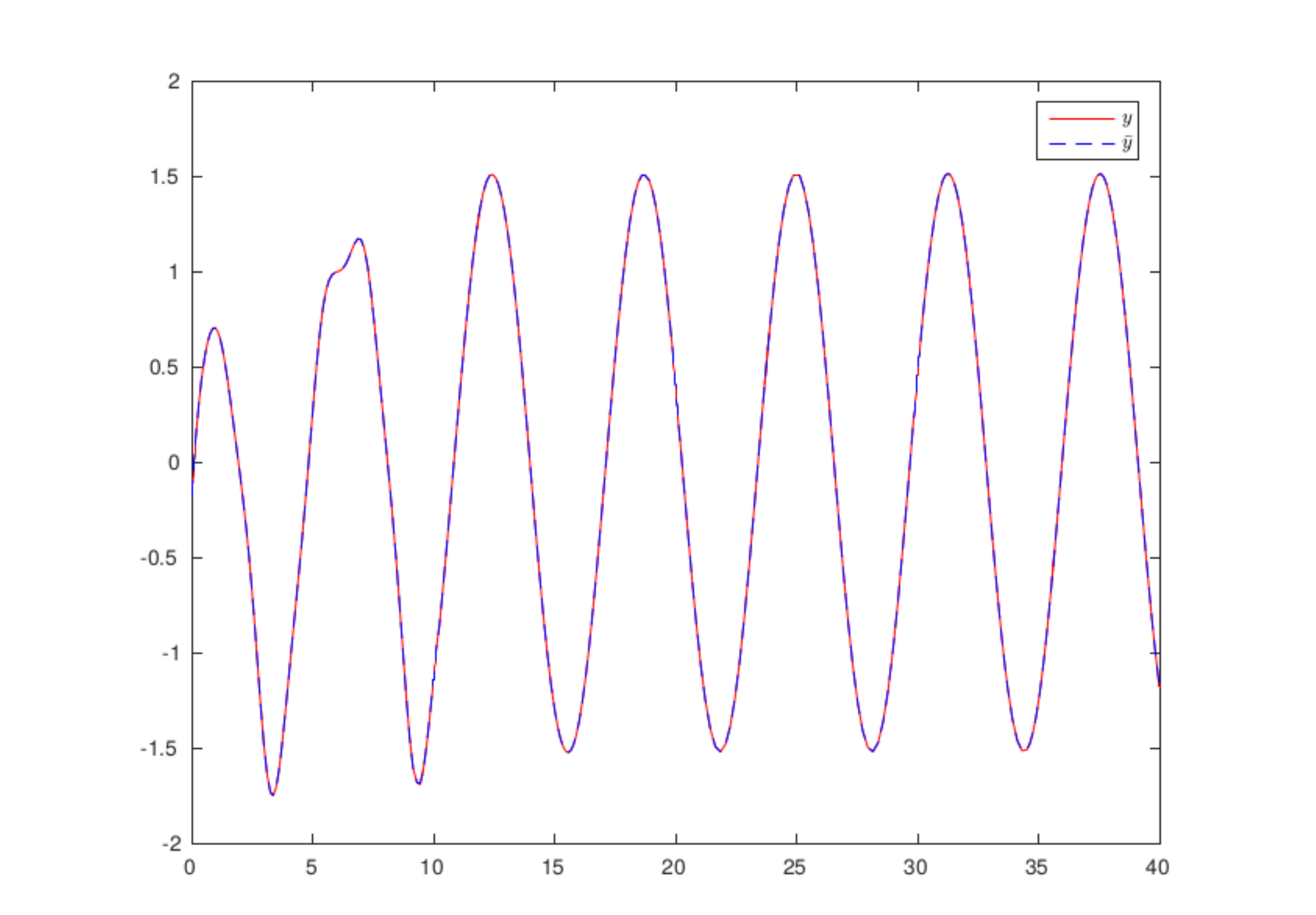}
	\caption{Comparisons of responses of $\Sigma$ (red) and $\bar{\Sigma}$ (green) for $n=200$, $r=200$, $m=4$. \label{fig3}}
  \end{figure}

  To demonstrate the scalability of the proposed approach, we tested it on a bilinear system
  $\Sigma$ of the form \eqref{eq:BSSform} with $m=4$ inputs and $n=200$ states and
  one output. In this case, $\gamma$ was chosen to be the language
  accepted by an NDFA of the form $\mathcal{A}=(S,Q,\{ \rightarrow_q \}_{q \in Q},F,s_0)$, where $S=\{1,2,\ldots,m\}=F$,  $s_0=1$ aand
 $\rightarrow_{0} =\{ (i,j) \mid i \le j \}$ and
 for all $q=1,\ldots,m$, $\rightarrow_{q} =\{ (i,q) \mid i \le q \}$.
  Hence, $L(\mathcal{A})=\gamma$.  More explicitly,
 $\gamma$ can be described as the set of all
  the sequences $w \in Q^{*}$ such that $w=v_1 \cdots v_m$, 
  $v_i \in \{0,i\}^{*}$, $i=1,\ldots,m$. 
  Applying Algorithm \ref{alg4} and Theorem \ref{theo:krylov1} to $\Sigma$ yielded a bilinear system 
  $\bar{\Sigma}=(1,4,r,\{\bar{A}\}_{q \in Q},\bar{C},\bar{x}_0)$ of order $r=9$.
  In this case, it took $4$ iterations for Algorithm \ref{alg4} to terminate, which is much smaller than the
  theoretical upper bound $200 \cdot 4$. 
  Due to lack of space, we do not present the matrices of $\Sigma$ and $\bar{\Sigma}$, the tex files with the matrices
  and the .mat files can be found among the supplementary material of this report. For each $i=0,\ldots; m$,  the matricec $A_i$ 
  are in the files with .tex and .mat extension called \texttt{NOLCOSBig\_exampeA}$i$. The matrix s $C$ is stored in the .tex and .mat file
  \texttt{NOLCOSBig\_exampeC}, and the initial state $x_0$ is in the file \texttt{NOLCOSBig\_exampex0}.
  In a similar manner, for each $i=0,\ldots,m$, $\bar{A}_i$ is stored in files called \texttt{NOLCOSBig\_exampeA}$i$r, and $\bar{C}$ and $\bar{x}_0$ are
  stored in the files named \texttt{NOLCOSBig\_exampeCr}, and 
\texttt{NOLCOSBig\_exampex0r}. 
\\
  By Theorem \ref{theo1}, $\bar{\Sigma}$ has the same output as $\Sigma$ on $[0,T]$,
  $T=50$, where for the input $u$ which satisfies 
  $u(t)=\cos(\pi t) e_i$ if $t \in [t_{i-1},t_{i})$ for all $i=1,\ldots,m$, where $t_{0}=0$ and $t_{i+1}=t_i+10$, $i=0,\ldots,m-1$.
  This complies with Theorem \ref{theo1}, as $u \in \mathcal{U}_{\gamma,T}$. The responses 
  $y=Y_{\Sigma,x_0}(u)$ and $\bar{y}=Y_{\bar{\Sigma},\bar{x}_0}(u)$ are shown on Figure \ref{fig3}, and 
  it can be seen that they are indeed the same. 
	
\begin{figure}[hbpt]
	\centering
  \includegraphics[scale=0.6]{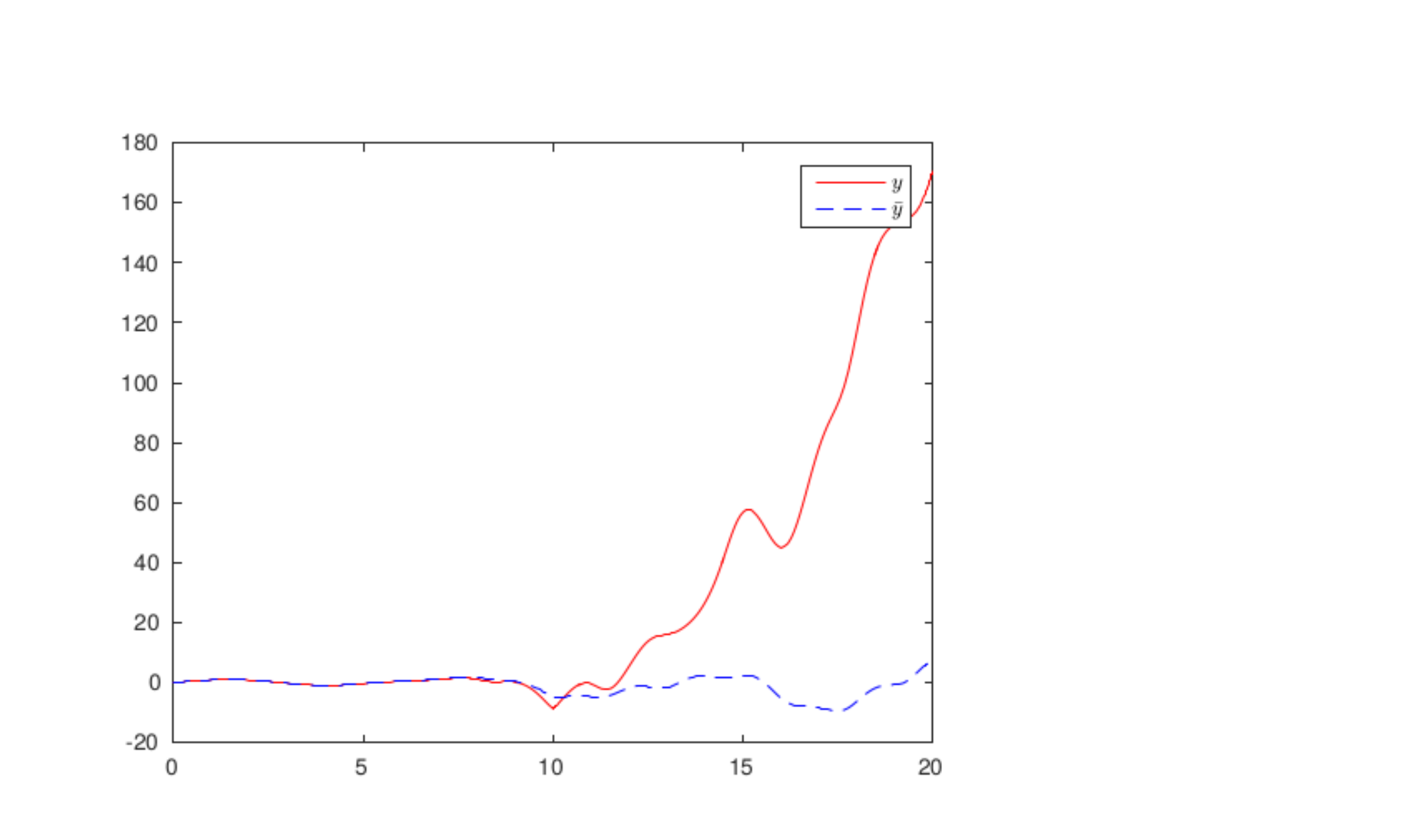}
	\caption{Comparisons of responses of $\Sigma$ (red) and $\bar{\Sigma}$ (green) for $n=200$, $r=200$, $m=4$. \label{fig4}} 
  \end{figure}

  We also applied Algorithm \ref{alg1} and Theorem \ref{theo:krylov1} to $\Sigma$ for $\gamma=\{v \in Q^{*} \mid |v| \le N\}$
  with $N=3$. In this case,  we obtained a reduced-order system $\hat{\Sigma}=(1,4,n_r,\{\hat{A}_q\}_{q \in Q},\hat{C},\hat{x}_0)$ of order $n_r=17$. 
  The matrices $\hat{A}_q$, $q \in Q$, $\hat{C}$ and the vector $\hat{x}_0$ can be found in the supplementary material, in the files named \texttt{NOLCOSBig\_exampeA$q$pr}, 
  \texttt{NOLCOSBig\_exampeCpr},\texttt{NOLCOSBig\_exampex0pr} respectively.
We evaluated the response $y$ and
  $\hat{y}$ of $\Sigma$ and $\bar{\Sigma}$ respectively for the following input:
  $u(t)=(1+\cos(\pi t))e_{q_i}+(\cos(t),\cos(2t),\ldots, \cos(mD))^T$ if $t \in [t_{i-1},t_i)$, $i=1,\ldots,k$, where $k=5$,
  $q_1=4,q_2=3, q_3=1, q_4=2,q_5=4$, $t_0=0$ and $t_{i+1}=t_i+10$, $i=0,\ldots,k-1$, see Figure \ref{fig4}.
  We observe that $y$ and $\hat{y}$ are not the same, but on $[0,T]$, $T=10$ they are close as predicted by Theorem \ref{theo2}.
  The code for this example can be found in the supplementary material, in the file \texttt{BilinearModelExample1.m}. 
  
  Finally, we evaluated the proposed method on the nonlinear circuit investigated in \cite{BreitenDamm,Bai2006406} with $N_o=150$ nonlinear resistors.
  Following these references, we applied Carleman bilinearization to the original model. As the results, the obtained bilinear system $\Sigma=(1,2,n,\{A_i\}_{i=0}^{1},C,x_0)$ 
   is of the order $n=N_o+N_o^2+1$. 
  The code for generating the model can be found in the supplementary material, in the file \texttt{NonlinearCircuitBilinModel.m}. In 
  order to obtain matrices $C$,$A_0,A_1,x_0$, the code in \texttt{NonlinearCircuitBilinModel.m} should be run. After the code has been run,
  the matrices $C$, $A_0,A_1,x_0$ will be stored in the files
\texttt{CNonlinRC.mat},  \texttt{ANonlinRC.mat}, \texttt{NNonlinRC.mat}, \texttt{x0NonlinRC.mat} respectively.

  We applied Algorithm \ref{alg5} and Theorem \ref{theo:krylov2} with $\alpha=\{\epsilon, 0,1,10,110,1110,1110\}$, the resulting reduced order bilinear
  system 
$\bar{\Sigma}=(1,2,6,\{\bar{A}_i\}_{i=0}^{1},\bar{C},\bar{x}_0)$ 
 (of order $6$). The matrices $\bar{C}$,$\bar{A}_0,\bar{A}_1,\bar{x}_0$ can be found in the supplementary material, in the files
\texttt{CNonlinRCnr.mat},  \texttt{ANonlinRCnr.mat}, \texttt{NNonlinRCnr.mat}, \texttt{x0NonlinRCnr.mat} respectively.
Algorithm \ref{alg5} took $5$ iterations to stop, which is much less than
the theoretical upper bound of $N_o|S|$, where $|S|$ is the number of states of the NDFA accepting $\alpha$.
Note that $N_o=150$ and $|S|=7$ in this case. 
We applied Theorem \ref{theo:krylov2}
  and Algorithm \ref{alg2} with $N=3$, which yielded a reduced order system
  $\hat{\Sigma}=(1,2,12,\{\hat{A}_i\}_{i=0}^{1},\hat{C}, \hat{x}_0)$ (of order 12). 
 The matrices $\hat{C}$,$\hat{A}_0,\hat{A}_1,\hat{x}_0$ can be found in the supplementary material, in the files
\texttt{CNonlinRCr.mat},  \texttt{ANonlinRCr.mat}, \texttt{NNonlinRCr.mat}, \texttt{x0NonlinRCr.mat}
 respectively.

\begin{figure}[hbpt]
 	\centering
  \includegraphics[scale=0.45]{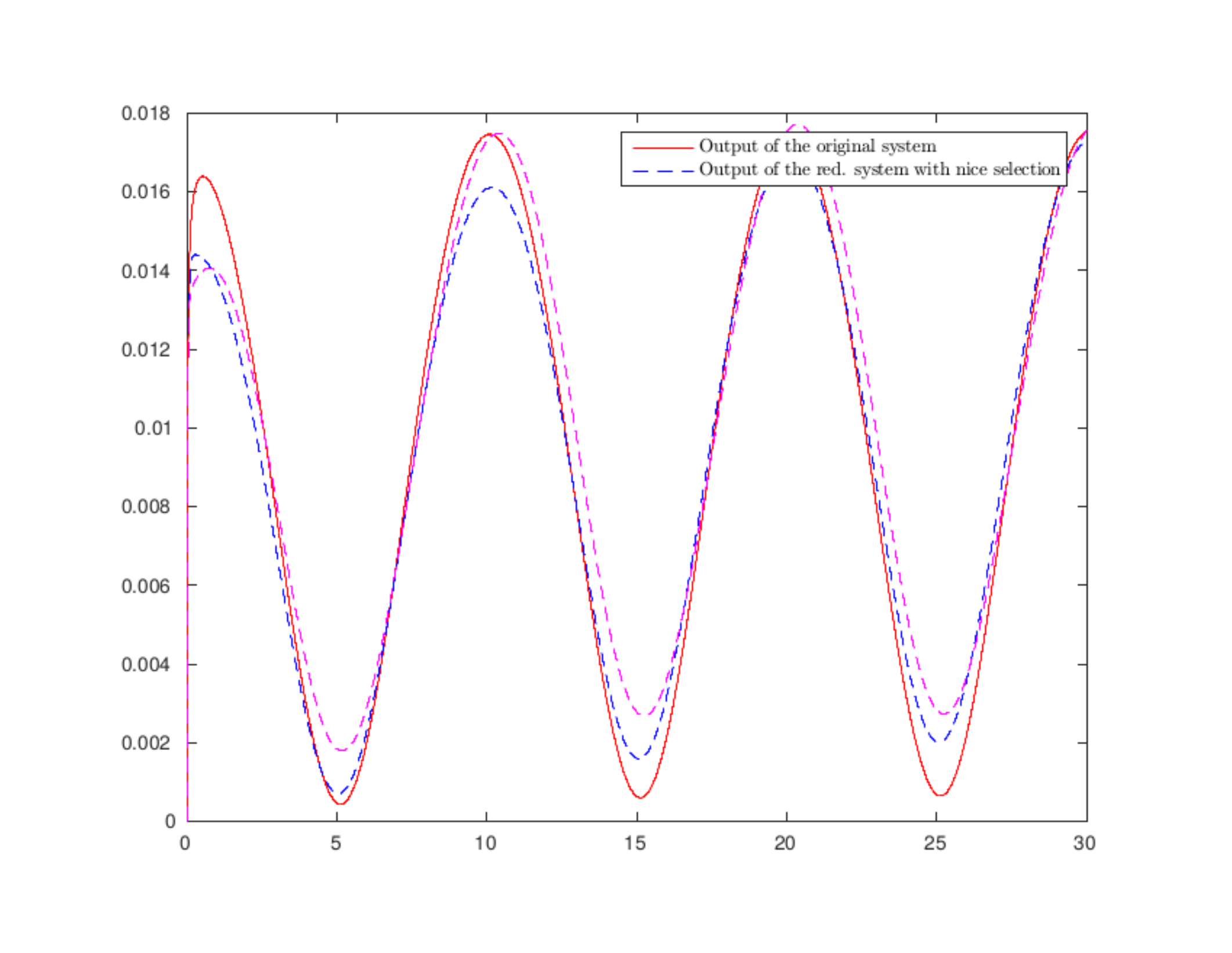}
	\caption{Comparisons of responses of the nonlinear RC and
  $\bar{y}$ and $\hat{y}$ of $\bar{\Sigma}$ and $\hat{\Sigma}$ \label{fig5}
}
  \end{figure}

We simulated the responses of the original nonlinear model 
  (before Carleman's bilinearization), and the responses
  $\bar{y}=Y_{\bar{\Sigma},\bar{x}_0},\hat{y}=Y_{\Sigma,\hat{x}_0}$ of $\bar{\Sigma}$ and $\hat{\Sigma}$ for the input 
 $u(t)=(\cos(2\pi t/10)+1)/2$, the result is shown on
Figure \ref{fig5}. The responses of $\bar{\Sigma}$ and $\hat{\Sigma}$ are
both reasonably close to the response of the original  nonlinear model, 
yet the order of $\bar{\Sigma}$ is much lower than that
of $\hat{\Sigma}$. This demonstrates that nice selections give additional 
flexibility to model reduction.
   The code for this example can be found in the supplementary material, in the file \texttt{TestNPartialNonlinCircuit.m}.

\section{Conclusion}
We have developed a method for model reduction of bilinear control systems leaning upon the concept of 
 the column nice selection and the row nice selection.  The resulting bilinear system has exactly the same output response as the original system for inputs consistent with a nice selection. For other inputs, the error between the time responses of the two systems decreases with the cardinality of the nice selection, provided the inputs and considered time horizon are short. Furthermore, we have provided algorithms for computing matrix representations of $\alpha$-unobservability and $\beta$-reachability spaces, which has been used for computing $\alpha$-partial and $\beta$-partial realizations of an input-output map. 
Future research will be directed towards a better understanding of the numerical issues involved, of error bounds for the reduced model, and of the relative advantage of the proposed method in comparison to \cite{Bai2006406,BilinearMomentMatching1,BilinearMomentMatching2,BilinearMomentMatching5,BreitenDamm,BennerSiam,FlagGugercin,Wang20121231}. 
\\
\textbf{Acknowledgement:} This work was partially supported by ESTIREZ project of Region Nord-Pas de Calais, France, and by the Innovation Fund Denmark, EDGE project (contract no. 11-116843).

\end{document}